\documentclass[fullpage,twoside]{article}
\usepackage{fancyhdr,graphicx}
\usepackage{amsfonts}
\usepackage{booktabs}
\usepackage{array}
\usepackage{amsmath}
\usepackage{graphicx}
\usepackage{amsmath,bm}
\def \u{{\mbox{\boldmath $u$}}}

\begin{document}
	\title{\vspace{-1cm} Interactive Fuzzy Goal Programming Based on Taylor Series to Solve Multiobjective Nonlinear Programming Problems with Interval Type 2 Fuzzy Numbers}
	\author{Hasan Dalman \\
		{\small\em Deparment of Mathematics Engineering, Yildiz Technical University, Esenler 34210, Turkey}\\
		{\small\em hsandalman@gmail.com}}
	\date{}
	\maketitle
	\maketitle \thispagestyle{fancy} \label{first}
	\vspace{-5mm} \pagestyle{plain}
	\newsavebox{\mytable}

\begin{abstract}
	%\vskip 0.15mm\noindent {\bf Abstract:}
 This paper presents an interactive fuzzy goal programming (FGP) approach for solving multiobjective nonlinear programming problems (MONLPP) with interval type 2 fuzzy numbers (IT2 FNs). The cost and time of the objective functions, the resources, and the requirements of each kind of resources are taken to be trapezoidal IT2 FNs. Here, the considered problem is first transformed into an equivalent crisp MONLPP, and then the transformed MONLPP is converted into an equivalent Multiobjective Linear Programming Problem (MOLPP). By using a procedure based on Taylor series, this problem is reduced into a single objective linear programming problem (LPP) which can be easily solved by Maple 18.02 optimization toolbox. Finally, the proposed solution procedure is illustrated by two numerical examples. 


\end{abstract}

\section{ Introduction}

Most of the real-life problems are frequently characterized by multiple and conflicting criteria. Such conditions are normally estimated by optimizing multiple objective functions. Besides, when modeling real-world problems, often the parameters are included inexact quantities due to varied unmanageable factors. In practical mathematical programming problems, a decision-maker generally encounters a situation of uncertainty as well as complexity, due to various unknown factors. Usually, it is required to optimize several nonlinear and conflicting objectives simultaneously. To address the uncertain parameters which result in such situations, different fuzzy numbers are employed.  Fuzzy quantities are very suitable for modeling these type conditions. 

The fuzzy set theory first developed by Zadeh [1]  has been used to decision-making problems with imprecise information. Bellman and Zadeh [2] introduce that a fuzzy decision making is defined as the fuzzy set of options, obtaining from the intersection of the goals or objectives and constraints. The concept of fuzzy programming was first introduced by Tanaka et al. [3] in the structure of the fuzzy decision of Bellman and Zadeh. Later, fuzzy programming approach to linear programming with many objectives was investigated by Zimmermann [4]. 

The simplest approach for solving the fuzzy linear programming problem is converted it into the corresponding crisp programming problem. Zimmermann [4] has developed a fuzzy programming approach to solve the crisp multi-objective linear programming problem. Some authors have transformed the fuzzy programming problem into the crisp problem by using the ranking function [5, 6] and then solved it by conventional methods. 

In many practical problems such as in industrial planning, financial and corporate planning, marketing and media selection, etc., there exist many fuzzy and nonlinear production, planning and scheduling problems. These problems cannot be expressed and solved by conventional techniques due to uncertain information. So, the investigation on modeling and optimization for nonlinear programming with interval type 2 fuzzy  numbers (IT2 FNs) are not only significant in the fuzzy programming theory but also have a great and wide advantage in the application of the real-world practical problems of conflicting nature.

Type-2 fuzzy sets are introduced by Zadeh et al. [7] as the extension of type-1 fuzzy sets. Type-2 fuzzy sets are characterized by two memberships to determine more degrees. Since type-2 fuzzy sets have the advantage of modeling uncertain systems more correctly compared with type-1 fuzzy sets. However, when the type-2 fuzzy sets are employed to solve the problems, the computational procedures are very complicated [8]. So interval type 2 (IT2) fuzzy sets are widely employed with some relative illustrations to decrease dimensions, which are highly useful for computation and theoretical studies [9]. IT2 fuzzy sets can be observed as a particular illustration of common type-2 fuzzy sets that all the values of secondary membership are equal to 1. Hence, it not only represents the uncertainty better than type-1 fuzzy sets but also reduces the complexity compared to type-2 fuzzy sets.

Mendel et al. presented some definitions and concepts of IT2 fuzzy sets in [8]. Mitchell [10] and Zeng and Li [11] suggested methods to describe the connection between IT2 fuzzy sets. To accomplish limitations in these methods, Wu and Mendel [12] suggested a method called vector similarity method to convert IT2 fuzzy sets into the word more effectively. Ondrej and Milos [13] used IT2 fuzzy sets to generate a fuzzy voter design for fault-tolerant systems. Shu and Liang [14] proposed a new method based on IT2  Fuzzy Logic Systems (FLSs) to investigate and evaluate the network lifetime for wireless sensor networks. Wu and Mendel [15] defined linguistic weighted average and used it to handle hierarchical multicriteria decision-making problems. Han and Mendel [16] IT2 FNs in deciding the logistic location, and the result has been demonstrated to be more comforting. Chen and Lee [17] proposed the definition of possibility degree of trapezoidal IT2 FNs and some arithmetic operations. Sinha et. al. [18] used IT2 FNs for modeling a multiobjective solid transportation problem. Li et. al. [19] investigated the problem of filter design for IT2 FLSs with D stability constraints based on a new performance index.  Up to now, IT2 FNs were used by many authors for decision-making problems [20, 21, 22, 23, 24, 25]. Due to their facility to handle with the high level of uncertainty, IT2 FLSs further performed in various real-world applications, containing intelligent control [26, 27], time series predictions [28, 29, 30], pattern recognition [31], image processing [32] and many others.

In this paper, an interactive fuzzy goal programming (FGP) approach based on Taylor series is presented to achieve the highest degree of membership function for multiobjective nonlinear programming (MONLPP) with trapezoidal IT2 FNs. The FGP approach first introduced  by Narasimhan [33] and then Hannan [34] presented different membership functions, i.e., piecewise linear membership functions into FGP model. Tiwari et al. [35] introduced a weighted additive model that incorporates each goal's weight into the objective function, where the weights reveal the relative importance of the fuzzy goals. Mohamed [36] discussed the relationship between goal programming and fuzzy programming where the highest degree of each of the membership goals is achieved by minimizing over deviation variables. Several methods are suggested to linearize the fractional and/or nonlinear functions in literature. In the case of a nonlinear programming, the most common methods are based on linearization procedures [37, 38, 39, 40].  

Here, all the parameters of MONLPP are considered trapezoidal IT2 FNs. To the best of my knowledge, no work has been studied on MONLPP with trapezoidal IT2 FNs under the nonlinear constraints. In order to convert the considered fuzzy model into its crisp equivalent, the expected value of trapezoidal IT2 FNs is first employed and then aspiration levels and the tolerance limits of the objective functions are determined by getting individual optimal solutions and thereby the feasible region for the problem is reconstructed by using the upper of lower limits of decision variables. After these operations, the nonlinear membership functions, which are associated with each nonlinear objective of the problem are constructed and then with the use of Taylor series approach around its maximal solution, each nonlinear membership function is converted into linear functions. In this way, this problem is reduced into a single objective linear programming problem and then an interactive solution procedure is presented to determine the optimal solution for MNLOPP with IT2 FNs. Finally, numerical examples are presented to demonstrate the feasibility of the suggested procedures. 

The paper is constructed as follows: Sect. 2 deals with some definitions and arithmetic operations on IT2 FNs. Section 3 deals with problem formulation and its solution procedure. In Sect. 4, numerical examples are given to illustrate the methodology. Finally, we concluded in Sect. 5.

\section{ Preliminaries }

\subsection{ Interval Type-2 Fuzzy Set}

\noindent \textbf{Definition 1 }(Mendel et. al. [8]) Let $\tilde{A}$ be a type-2 fuzzy set, then $\tilde{A}$is defined as

\noindent $\tilde{A}=\left\{\left. \left(x,\mu \right),\mu _{\tilde{A}} \left(x,\mu \right)\right|\forall x\in X,\forall \mu \in {\rm {\rm F}}\subseteq \left[0,1\right],0\le \mu _{\tilde{A}} \left(x,\mu \right)\le 1\right\},$where $X$is the universe of discourse and $\mu _{\tilde{A}} \left(x,\mu \right)$denotes the membership function of $\tilde{A}.$ $\tilde{A}$ can be defined as $\tilde{A}=\int _{x\in X}\int _{\mu \in {\rm {\rm F}}}{\raise0.7ex\hbox{$ \mu _{\tilde{A}} \left(x,\mu \right) $}\!\mathord{\left/ {\vphantom {\mu _{\tilde{A}} \left(x,\mu \right) \left(x,\mu \right)}} \right. \kern-\nulldelimiterspace}\!\lower0.7ex\hbox{$ \left(x,\mu \right) $}}   \mathop{}\nolimits_{,}^{} \mu \in {\rm {\rm F}}\subseteq \left[0,1\right].$

\noindent \textbf{Definition 2 }(Mendel et. al. [8]) If all $\mu _{\tilde{A}} \left(x,\mu \right)=1,$then $\tilde{A}$called an IT2 fuzzy set i.e. $\tilde{A}=\int _{x\in X}\int _{\mu \in {\rm {\rm F}}}{\raise0.7ex\hbox{$ 1 $}\!\mathord{\left/ {\vphantom {1 \left(x,\mu \right)}} \right. \kern-\nulldelimiterspace}\!\lower0.7ex\hbox{$ \left(x,\mu \right) $}}  , \mathop{}\nolimits_{}^{} \mu \in {\rm {\rm F}}\subseteq \left[0,1\right].$

Uncertainty in the first memberships of a type-2 fuzzy set $\tilde{A}$consists of a bounded region that we call the footprint of uncertainty. It is the union of all first memberships. 

The footprint of uncertainty is characterized by the upper membership function and the lower membership function, and are denoted by $\bar{\mu }_{\tilde{A}} $and $\underline{\mu }_{\tilde{A}} $(Mendel et. al. [8]).

\textbf{Definition 3 }An IT2 FN is called a trapezoidal IT2 FN where the upper membership function and the lower membership function are both trapezoidal fuzzy numbers, i.e., 
\begin{equation} \label{GrindEQ__1_} 
A=\left(\bar{A},\underline{A}\right)=\left(\bar{a}_{1}^{} ,\bar{a}_{2}^{} ,\bar{a}_{3}^{} ,\bar{a}_{4}^{} ;H_{1} \left(\bar{A}\right),H_{2} \left(\bar{A}\right)\right), 
\end{equation} 
\[\left(\underline{a}_{1}^{} ,\underline{a}_{2}^{} ,\underline{a}_{3}^{} ,\underline{a}_{4}^{} ;H_{1} \left(\underline{A}\right),H_{2} \left(\underline{A}\right)\right),\] 
where $H_{j} \left(\underline{A}\right)$ and $H_{j} \left(\bar{A}\right)$ denote membership values of the corresponding elements $\underline{a}_{j+1}^{} $ and $\bar{a}_{j+1}^{} $respectively.

\noindent

\subsection{ The Arithmetic Operations of Interval Type-2 Fuzzy Set }

Suppose $A_{1} $ and $A_{2} $ are two trapezoidal IT2 FNs, then the following procedures are satisfied (Li et. al [19]):
\begin{equation} \label{GrindEQ__2_} 
A_{1} =\left(\bar{A}_{1} ,\underline{A}_{1} \right)=\left(\bar{a}_{11}^{} ,\bar{a}_{12}^{} ,\bar{a}_{13}^{} ,\bar{a}_{14}^{} ;H_{1} \left(\bar{A}\right),H_{2} \left(\bar{A}\right)\right), 
\end{equation} 
\[\left(\underline{a}_{11}^{} ,\underline{a}_{12}^{} ,\underline{a}_{13}^{} ,\underline{a}_{14}^{} ;H_{1} \left(\underline{A}\right),H_{2} \left(\underline{A}\right)\right),\] 
\begin{equation} \label{GrindEQ__3_} 
A_{2} =\left(\bar{A}_{2} ,\underline{A}_{2} \right)=\left(\bar{a}_{21}^{} ,\bar{a}_{22}^{} ,\bar{a}_{23}^{} ,\bar{a}_{24}^{} ;H_{1} \left(\bar{A}\right),H_{2} \left(\bar{A}\right)\right), 
\end{equation} 
\[\left(\underline{a}_{21}^{} ,\underline{a}_{22}^{} ,\underline{a}_{23}^{} ,\underline{a}_{24}^{} ;H_{1} \left(\underline{A}\right),H_{2} \left(\underline{A}\right)\right),\] 
\begin{equation} \label{GrindEQ__4_} 
\begin{array}{rcl} {A_{1} +A_{2} } & {=} & {\left(\bar{A}_{1} ,\underline{A}_{1} \right)+\left(\bar{A}_{2} ,\underline{A}_{2} \right)} \\ {} & {=} & {\left(\begin{array}{l} {\bar{a}_{11}^{} +\bar{a}_{21}^{} ,\bar{a}_{12}^{} +\bar{a}_{22}^{} ,\bar{a}_{13}^{} +\bar{a}_{23}^{} ,\bar{a}_{14}^{} +\bar{a}_{24}^{} ;} \\ {\min \left(H_{1} \left(\bar{A}_{1} \right),H_{1} \left(\bar{A}_{2} \right)\right),\min \left(H_{2} \left(\bar{A}_{1} \right),H_{2} \left(\bar{A}_{2} \right)\right)} \\ {\underline{a}_{11}^{} +\underline{a}_{21}^{} ,\underline{a}_{12}^{} +\underline{a}_{21}^{} ,\underline{a}_{13}^{} +\underline{a}_{23}^{} ,\underline{a}_{14}^{} +\underline{a}_{24}^{} ;} \\ {\min \left(H_{1} \left(\underline{A}_{1} \right),H_{1} \left(\underline{A}_{2} \right)\right),\min \left(H_{2} \left(\underline{A}_{1} \right),H_{2} \left(\underline{A}_{2} \right)\right),} \end{array}\right)} \end{array} 
\end{equation} 
\begin{equation} \label{GrindEQ__5_} 
\begin{array}{rcl} {A_{1} -A_{2} } & {=} & {\left(\bar{A}_{1} ,\underline{A}_{1} \right)-\left(\bar{A}_{2} ,\underline{A}_{2} \right)} \\ {} & {=} & {\left(\begin{array}{l} {\bar{a}_{11}^{} -\bar{a}_{21}^{} ,\bar{a}_{12}^{} -\bar{a}_{22}^{} ,\bar{a}_{13}^{} -\bar{a}_{23}^{} ,\bar{a}_{14}^{} -\bar{a}_{24}^{} ;} \\ {\min \left(H_{1} \left(\bar{A}_{1} \right),H_{1} \left(\bar{A}_{2} \right)\right),\min \left(H_{2} \left(\bar{A}_{1} \right),H_{2} \left(\bar{A}_{2} \right)\right)} \\ {\underline{a}_{11}^{} -\underline{a}_{21}^{} ,\underline{a}_{12}^{} -\underline{a}_{21}^{} ,\underline{a}_{13}^{} -\underline{a}_{23}^{} ,\underline{a}_{14}^{} -\underline{a}_{24}^{} ;} \\ {\min \left(H_{1} \left(\underline{A}_{1} \right),H_{1} \left(\underline{A}_{2} \right)\right),\min \left(H_{2} \left(\underline{A}_{1} \right),H_{2} \left(\underline{A}_{2} \right)\right),} \end{array}\right)} \end{array} 
\end{equation} 
\[A_{1} \times A_{2} =\left(\bar{A}_{1} ,\underline{A}_{1} \right)\times \left(\bar{A}_{2} ,\underline{A}_{2} \right)\] 
\begin{equation} \label{GrindEQ__6_} 
=\left(\begin{array}{l} {\bar{a}_{11}^{} \times \bar{a}_{21}^{} ,\bar{a}_{12}^{} \times \bar{a}_{22}^{} ,\bar{a}_{13}^{} \times \bar{a}_{23}^{} ,\bar{a}_{14}^{} \times \bar{a}_{24}^{} ;} \\ {\min \left(H_{1} \left(\bar{A}_{1} \right),H_{1} \left(\bar{A}_{2} \right)\right),\min \left(H_{2} \left(\bar{A}_{1} \right),H_{2} \left(\bar{A}_{2} \right)\right)} \\ {\underline{a}_{11}^{} \times \underline{a}_{21}^{} ,\underline{a}_{12}^{} \times \underline{a}_{21}^{} ,\underline{a}_{13}^{} \times \underline{a}_{23}^{} ,\underline{a}_{14}^{} \times \underline{a}_{24}^{} ;} \\ {\min \left(H_{1} \left(\underline{A}_{1} \right),H_{1} \left(\underline{A}_{2} \right)\right),\min \left(H_{2} \left(\underline{A}_{1} \right),H_{2} \left(\underline{A}_{2} \right)\right),} \end{array}\right) 
\end{equation} 
\begin{equation} \label{GrindEQ__7_} 
\begin{array}{rcl} {kA_{1} =\left(\bar{A}_{1} ,\underline{A}_{1} \right)=\left(k\bar{a}_{11}^{} ,k\bar{a}_{12}^{} ,k\bar{a}_{13}^{} ,k\bar{a}_{14}^{} ;H_{1} \left(\bar{A}_{2} \right),H_{2} \left(\bar{A}_{2} \right)\right)} & {,} & {} \\ {\left(k\underline{a}_{11}^{} ,k\underline{a}_{12}^{} ,k\underline{a}_{13}^{} ,k\underline{a}_{14}^{} ;H_{1} \left(\underline{A}_{2} \right),H_{2} \left(\underline{A}_{1} \right)\right)} & {,} & {} \end{array} 
\end{equation} 
\begin{equation} \label{GrindEQ__8_} 
\begin{array}{rcl} {\frac{1}{k} A_{1} =\left(\bar{A}_{1} ,\underline{A}_{1} \right)=\left(\frac{1}{k} \bar{a}_{11}^{} ,\frac{1}{k} \bar{a}_{12}^{} ,\frac{1}{k} \bar{a}_{13}^{} ,\frac{1}{k} \bar{a}_{14}^{} ;H_{1} \left(\bar{A}_{2} \right),H_{2} \left(\bar{A}_{2} \right)\right)} & {,} & {} \\ {\left(\frac{1}{k} \underline{a}_{11}^{} ,\frac{1}{k} \underline{a}_{12}^{} ,\frac{1}{k} \underline{a}_{13}^{} ,\frac{1}{k} \underline{a}_{14}^{} ;H_{1} \left(\underline{A}_{2} \right),H_{2} \left(\underline{A}_{1} \right)\right)} & {,} & {} \end{array} 
\end{equation}

\subsection{ Defuzzification of Trapezoidal Interval Type-2 Fuzzy Numbers}

Let us consider a trapezoidal IT2 FN$A$characterized by Equation \eqref{GrindEQ__1_}. The expected value of $A$ is determined as follows (Hu et. al [20]):
\begin{equation} \label{GrindEQ__9_} 
E\left(A\right)=\frac{1}{2} \left(\frac{1}{4} \sum _{i=1}^{4}\left(\underline{a}_{i} +\bar{a}_{i} \right) \right)\times \frac{1}{4} \left(\sum _{i=1}^{2}\left(H_{i} \left(\underline{A}_{i} \right)+H_{i} \left(\bar{A}_{i} \right)\right) \right) 
\end{equation} 

Assuming that $A_{1} $ and $A_{2} $ are two trapezoidal IT2 FNs, then we get $A_{1} >A_{2} $if and only if $E\left(A_{1} \right)>E\left(A_{2} \right).$

\noindent When$\bar{a}_{i}^{} =\underline{a}_{i}^{} ,$$(i=1,2,3,4)$and $H_{1} \left(\underline{A}\right)=H_{2} \left(\underline{A}\right)=H_{1} \left(\bar{A}\right)=H_{2} \left(\bar{A}\right)$ the trapezoidal IT2 FN reduces to trapezoidal fuzzy number, just as $\tilde{A}=\left(\underline{a}_{1}^{} ,\underline{a}_{2}^{} ,\underline{a}_{3}^{} ,\underline{a}_{4}^{} \right).$The expected value of $\tilde{A}$is 
\[E\left(\tilde{A}\right)=\left({\raise0.7ex\hbox{$ \underline{a}_{1}^{} +\underline{a}_{2}^{} +\underline{a}_{3}^{} +\underline{a}_{4}^{}  $}\!\mathord{\left/ {\vphantom {\underline{a}_{1}^{} +\underline{a}_{2}^{} +\underline{a}_{3}^{} +\underline{a}_{4}^{}  4}} \right. \kern-\nulldelimiterspace}\!\lower0.7ex\hbox{$ 4 $}} \right).\]

\section{ Problem Formulation}

A conventional MONLPP is formulated as:
\begin{equation} \label{GrindEQ__10_} 
\begin{array}{l} {Opt\, f_{k} \left(x\right),\, \, k=1,2,...,l_{} ,} \\ {s.t.\left\{\begin{array}{l} {g_{j} \left(x\right)\le b_{j} ,\, \, \, \, \, j=1,2,...,m_{1} .} \\ {g_{j} \left(x\right)\ge b_{j} ,\, \, \, \, \, j=m_{1} +1,m_{1} +2,...,m_{2} ,} \\ {g_{j} \left(x\right)=b_{j} ,\, \, \, \, j=m_{2} +1,m_{2} +2,...,m_{} ,} \\ {x\ge 0} \end{array}\right. } \end{array} 
\end{equation} 
where ``$Opt$'' denotes minimization and maximization; $f_{k} \left(x\right),\, \, \, k=1,2,...,l$ are multiple and nonlinear objectives to be optimized;$g_{j} \left(x\right),\, \, j=1,2,...,m$ are real-valued nonlinear constraints.$x=\left(x_{1}^{} ,x_{2}^{} ,...,x_{n}^{} \right)\in X$is $n-$ dimensionel decision vector.

In real-world decision-making problems such as in production, planning, scheduling, etc. the present quantity of resources as well as the production quantity or the demand quantity or the target over a period might be imprecise and possess various types of fuzziness due to many factors such as market price, existence of men power, perception with the operators, weather, rain, transportation, traffic, etc. Also, the objectives characterized by the decision-maker may be ill-defined due to estimated parameters. Thus, IT2 FNs appears to be more practical in such conditions. In formulating such problems, the detailed concepts and notations are given in papers [41, 42].

Assuming that the objective functions $\tilde{f}_{k} \left(x\right)$ and the resource constraint functions $g_{j} \left(x\right)$ are nonlinear with estimated coefficient parameters which are in terms of trapezoidal IT2 FNs. The MONLPP with trapezoidal IT2 FNs can be formulated as:
\begin{equation} \label{GrindEQ__11_} 
\begin{array}{l} {Opt.\, \, \tilde{f}_{k} \left(x\right),\, \, k=1,2,...,l,} \\ {s.t.\left\{\begin{array}{l} {\tilde{g}_{j} \left(x\right)\le \tilde{b}_{j} ,\, \, \, \, \, j=1,2,...,m_{1} .} \\ {\tilde{g}_{j} \left(x\right)\ge \tilde{b}_{j} ,\, \, \, \, \, j=m_{1} +1,m_{1} +2,...,m_{2} ,} \\ {\tilde{g}_{j} \left(x\right)=\tilde{b}_{j} ,\, \, \, \, \, j=m_{2} +1,m_{2} +2,...,m_{} ,} \\ {x\ge 0} \end{array}\right. } \end{array} 
\end{equation} 
where $\tilde{f}_{k} \left(x\right)=\sum _{i=1}^{l_{i} }\tilde{c}_{ik}^{}  \prod _{l=1}^{n}x_{l}^{\alpha _{l} }  ,k=1,2,...,l;$ $\tilde{g}_{j} \left(x\right)=\sum _{r=1}^{r_{j} }\tilde{a}_{rj}^{}  \prod _{l=1}^{n}x_{l}^{\beta _{l} }  ,j=1,2,...,m;$$\tilde{c}_{ik}^{} ,\, k=1,2,...,l;$$\tilde{a}_{rj}^{} ,\, \, j=1,2,...,m;$$\tilde{b}_{j} ,\, \, \, \, \, j=1,2,...,m_{1} $ are considered to be IT2 FNs. $x$ is $n-$ dimensionel decision variable vector $x_{}^{} =\left(x_{1}^{} ,x_{2}^{} ,...,x_{n}^{} \right).$ Here, $\tilde{c}_{ik}^{} ,\, k=1,2,...,l$is the estimated coeeficient parameters of the objective functions; $\tilde{a}_{rj}^{} ,\, \, j=1,2,...,m$is the requirements of resources;$\tilde{b}_{j} ,\, \, \, \, \, j=1,2,...,m_{1} $ is the available quantity of resources,respectively.

\noindent Therefore, MONLPP with trapezoidal IT2 FNs \eqref{GrindEQ__11_} can be genarally formulated as follows:
\begin{equation} \label{GrindEQ__12_} 
\begin{array}{l} {Opt.\, \, \tilde{f}_{k} \left(x\right)=\sum _{i=1}^{l_{i} }\tilde{c}_{ik}^{}  \prod _{l=1}^{n}x_{l}^{\alpha _{l} }  ,\, \, k=l_{1} +1,l_{1} +2,...,l,} \\ {s.t.\left\{\begin{array}{l} {\sum _{r=1}^{r_{j} }\tilde{a}_{rj}^{}  \prod _{l=1}^{n}x_{l}^{\beta _{l} }  \le \tilde{b}_{j} ,\, \, j=1,2,...,m_{1} .} \\ {\sum _{r=1}^{r_{j} }\tilde{a}_{rj}^{}  \prod _{l=1}^{n}x_{l}^{\beta _{l} }  \ge \tilde{b}_{j} ,\mathop{}\nolimits_{} j=m_{1} +1,m_{1} +2,...,m_{2} ,} \\ {\sum _{r=1}^{r_{j} }\tilde{a}_{rj}^{}  \prod _{l=1}^{n}x_{l}^{\beta _{l} }  =\tilde{b}_{j} \mathop{}\nolimits_{} j=m_{2} +1,m_{2} +2,...,m_{} ,} \\ {x_{l} \ge 0,l=1,2,...,n} \end{array}\right. } \end{array} 
\end{equation} 
By using the expected value function as defined in \eqref{GrindEQ__9_},  problem \eqref{GrindEQ__12_} is transformed into an equivalent crisp MONLPP as:
\begin{equation} \label{GrindEQ__13_} 
\begin{array}{l} {Opt.\, \, f_{k} \left(x\right)\cong \sum _{i=1}^{l_{i} }c_{ik}^{}  \prod _{l=1}^{n}x_{l}^{\alpha _{l} }  ,k=l_{1} +1,l_{1} +2,...,l,} \\ {s.t.\left\{\begin{array}{l} {g_{j} \left(x\right)\cong \sum _{r=1}^{r_{j} }a_{rj}^{}  \prod _{l=1}^{n}x_{l}^{\beta _{l} }  \le b_{j} ,\mathop{}\nolimits_{} j=1,2,...,m_{1} .} \\ {g_{j} \left(x\right)\cong \sum _{r=1}^{r_{j} }a_{rj}^{}  \prod _{l=1}^{n}x_{l}^{\beta _{l} }  \ge b_{j} ,\mathop{}\nolimits_{} j=m_{1} +1,m_{1} +2,...,m_{2} ,} \\ {g_{j} \left(x\right)\cong \sum _{r=1}^{r_{j} }a_{rj}^{}  \prod _{l=1}^{n}x_{l}^{\beta _{l} }  =b_{j} \mathop{}\nolimits_{} j=m_{2} +1,m_{2} +2,...,m_{} ,} \\ {x_{l} \ge 0,l=1,2,...,n} \end{array}\right. } \end{array} 
\end{equation} 
where the expected values of $\tilde{c}_{ik}^{} $, $\tilde{a}_{rj}^{} $ and $\tilde{b}_{j} $are $c_{ik}^{} ,\, \, k=1,2,...,l,\, \, a_{rj}^{} ,\, \, \, j=1,2,...,m$ and $b_{j} ,\, \, j=1,2,...,m,$respectively.

\subsection{ Construction of Fuzzy Multiobjective Nonlinear Goal Programming}

\noindent In a multiobjective programming, if an imprecise aspiration level is injected to each of the objectives, then these fuzzy objectives are expressed as fuzzy goals. Let $s_{k} $ be the aspiration level assigned to the $k^{th} $objective $f_{k} \left(x\right).$Then the fuzzy goals are $f_{k} \left(x\right)\lower3pt\hbox{\rlap{$\scriptscriptstyle\sim$}}\succ s_{k} $ for the maximization objective in \eqref{GrindEQ__13_} and $f_{k} \left(x\right)\lower3pt\hbox{\rlap{$\scriptscriptstyle\sim$}}\prec s_{k} $for the minimization objective of \eqref{GrindEQ__13_} where $\lower3pt\hbox{\rlap{$\scriptscriptstyle\sim$}}\succ $ and $\lower3pt\hbox{\rlap{$\scriptscriptstyle\sim$}}\prec $ represent the fuzzified inequalities.

\noindent Therefore, the fuzzy multiobjective goal programming problem can be formulated as follows:
\begin{equation} \label{GrindEQ__14_} 
\begin{array}{l} {f_{k} \left(x\right)\lower3pt\hbox{\rlap{$\scriptscriptstyle\sim$}}\succ s_{k} ,k=1,2,...,l_{1} ,} \\ {f_{k} \left(x\right)\lower3pt\hbox{\rlap{$\scriptscriptstyle\sim$}}\prec s_{k} ,k=l_{1} +1,l_{1} +2,...,l,} \\ {s.t.\left\{\begin{array}{l} {g_{j} \left(x\right)\cong \sum _{r=1}^{r_{j} }a_{rj}^{}  \prod _{l=1}^{n}x_{l}^{\beta _{l} }  \le b_{j} ,\mathop{}\nolimits_{} j=1,2,...,m_{1} .} \\ {g_{j} \left(x\right)\cong \sum _{r=1}^{r_{j} }a_{rj}^{}  \prod _{l=1}^{n}x_{l}^{\beta _{l} }  \ge b_{j} ,\mathop{}\nolimits_{} j=m_{1} +1,m_{1} +2,...,m_{2} ,} \\ {g_{j} \left(x\right)\cong \sum _{r=1}^{r_{j} }a_{rj}^{}  \prod _{l=1}^{n}x_{l}^{\beta _{l} }  =b_{j} \mathop{}\nolimits_{} j=m_{2} +1,m_{2} +2,...,m_{} ,} \\ {x_{l} \ge 0,l=1,2,...,n} \end{array}\right. } \end{array} 
\end{equation} 
Now, consider the $k^{th} $ fuzzy goal $f_{k} \left(x\right)\lower3pt\hbox{\rlap{$\scriptscriptstyle\sim$}}\succ s_{k} .$ Its membership function can be defined as follows:
\begin{equation} \label{GrindEQ__15_} 
\mu _{k} \left(f_{k} \left(x\right)\right)\cong \left\{\begin{array}{l} {1\, \, \, \, \, \, \, \, \, \, \, \, \, \, \, \, \, \, \, \, \, \, \, \, \, \, \, \, f_{k} \left(x\right)\ge s_{k} } \\ {\frac{f_{k} \left(x\right)-L_{k} }{\left(s_{k} -L_{k} \right)} ,\, \, \, \, \, L_{k} \le f_{k} \left(x\right)\le s_{k} \, \, \, \, \, \, \, } \\ {0\, \, \, \, \, \, \, \, \, \, \, \, \, \, \, \, \, \, \, \, \, \, \, \, \, \, \, \, L_{k} \ge f_{k} \left(x\right)} \end{array}\right.  
\end{equation} 
where $L_{k} $is the lower tolerance limit for the $k^{th} $fuzzy goal  and $\left(s_{k} -L_{k} \right)$is the tolerant interval which is subjectively selected, respectively. They are determined as follows:

$s_{k} =Max\left\{f_{k} \left(x\right),x\in X\right\}$and $L_{k} =Min\left\{f_{k} \left(x\right),x\in X\right\},$$k=1,2,...,l.$                      \eqref{GrindEQ__16_}

\noindent Similarly, consider the $k^{th} $ fuzzy goal of $f_{k} \left(x\right)\lower3pt\hbox{\rlap{$\scriptscriptstyle\sim$}}\prec s_{k} $. Its membership function can be defined as follows:
\begin{equation} \label{GrindEQ__17_} 
\mu _{k} \left(f_{k} \left(x\right)\right)=\left\{\begin{array}{l} {1\, \, \, \, \, \, \, \, \, \, \, \, \, \, \, \, \, \, \, \, \, \, \, \, \, \, \, \, f_{k} \left(x\right)\le s_{k} } \\ {\frac{U_{k} -f_{k} \left(x\right)}{U_{k} -s_{k} } ,\, \, \, \, s_{k} \le f_{k} \left(x\right)\le U_{k} } \\ {0\, \, \, \, \, \, \, \, \, \, \, \, \, \, \, \, \, \, \, \, \, \, \, \, \, \, \, U_{k} \le f_{k} \left(x\right)} \end{array}\right. \begin{array}{l} {} \\ {} \\ {} \end{array} 
\end{equation} 
where $U_{k} $is the upper tolerance limit for the $k^{th} $ fuzzy goal  and $\left(U_{k} -s_{k} \right)$ the tolerant interval which is subjectively selected, respectively. They are determined as follows:

$U_{k} =Max\left\{f_{k} \left(x\right),x\in X\right\}$and $s_{k} =Min\left\{f_{k} \left(x\right),x\in X\right\},$$k=1,2,...,l.$                      \eqref{GrindEQ__18_}

\noindent By using the max-min form introduced by Zadeh with the membership function as defined in \eqref{GrindEQ__15_} and \eqref{GrindEQ__17_}, a crisp nonlinear programming problem can be formulated as follows:
\begin{equation} \label{GrindEQ__19_} 
\begin{array}{l} {Max\, \, \lambda } \\ {\left\{\begin{array}{l} {\, \lambda \le \frac{f_{k} \left(x\right)-L_{k} }{\left(s_{k} -L_{k} \right)} ,\, \, k=1,2,..,l,} \\ {\lambda \le \frac{U_{k} -f_{k} \left(x\right)}{U_{k} -s_{k} } ,\, \, k=1,2,..,l,} \\ {g_{j} \left(x\right)\le b_{j} ,\, j=1,2,...,m_{1} .} \\ {g_{j} \left(x\right)\ge b_{j} ,\, j=m_{1} +1,m_{1} +2,...,m_{2} ,} \\ {g_{j} \left(x\right)=b_{j} ,\, j=m_{2} +1,m_{2} +2,...,m_{} ,} \\ {\lambda \le 1,\, \, \, \lambda \ge 0\, ,\, \, x_{l}^{} \ge 0,\, \, \, l=1,2,...,n.} \end{array}\right. } \end{array} 
\end{equation} 
Similarly, by using the min- max form, a crisp nonlinear programming problem can be formulated as follows:
\begin{equation} \label{GrindEQ__20_} 
\begin{array}{l} {Min\, \, \, 1-\, \lambda } \\ {\left\{\begin{array}{l} {\, 1-\lambda \ge 1-\frac{f_{k} \left(x\right)-L_{k} }{\left(s_{k} -L_{k} \right)} ,\, \, k=1,2,..,l,} \\ {1-\lambda \ge 1-\frac{U_{k} -f_{k} \left(x\right)}{U_{k} -s_{k} } ,\, \, k=1,2,..,l,} \\ {g_{j} \left(x\right)\le b_{j} ,\, j=1,2,...,m_{1} .} \\ {g_{j} \left(x\right)\ge b_{j} ,\, j=m_{1} +1,m_{1} +2,...,m_{2} ,} \\ {g_{j} \left(x\right)=b_{j} ,\, j=m_{2} +1,m_{2} +2,...,m_{} ,} \\ {\left(1-\lambda \right)\le 1,\, \, \, \left(1-\lambda \right)\ge 0,\, \, x_{l}^{} \ge 0,\, \, \, l=1,2,...,n.\, } \end{array}\right. } \end{array} 
\end{equation}

\subsection{ Linearization Nonlinear Membership and Constraint Functions Using the Taylor Series}

Several methods are implemented to linearize the fractional and/or nonlinear functions in literature [37, 38, 39, 40]. In this section, problem \eqref{GrindEQ__13_} will transform into an equivalent multiobjective linear programming problem (MOLPP). 

Note that the feasible region for a programming problem is the whole set of alternatives for the decision variables over which the objective function is to be optimized. Therefore, it can be illustrated with the limits of decision variables, and thereby the nonlinear constrained region can be easily converted to the linear inequalities.  

\noindent The suggested solution procedure can be continued as follows:

\begin{enumerate}
\item  Construct problem \eqref{GrindEQ__13_}

\item  Solve problem \eqref{GrindEQ__13_} as a single objective nonlinear programming problem, taking each time only one objective as objective function and ignoring all others.

\item  Compute the value of each objective function at each solution and then define the feasible region by the limits of decision variables.

\item  Determine $\tilde{x}_{l}^{*} =\left(\tilde{x}_{1}^{*} ,\tilde{x}_{2}^{*} ,...,\tilde{x}_{n}^{*} \right)$ which is the solution that is employed to maximize the $k^{th} $ nonlinear membership function $\mu _{k} \left(f_{k} \left(x\right)\right)$  associated with $k^{th} $ nonlinear objective $f_{k} \left(x\right).$ 

\item  Then, transform nonlinear membership functions by using Taylor series approach around the solution $\tilde{x}_{l}^{*} =\left(\tilde{x}_{1}^{*} ,\tilde{x}_{2}^{*} ,...,\tilde{x}_{n}^{*} \right)$  as follows:
\begin{equation} \label{GrindEQ__21_} 
\tilde{\mu }_{k} \left(f_{k} \left(x\right)\right)_{k=1,2,...,l.} \cong \left[\left. \frac{\mu _{k} \left(f_{k} \left(\tilde{x}_{l}^{*} \right)\right)}{\partial x_{1} } \right|_{\tilde{x}_{l} ^{*} } \left(x_{1} -\tilde{x}_{1}^{*} \right)+\left. \frac{\mu _{k} \left(f_{k} \left(\tilde{x}_{l}^{*} \right)\right)}{\partial x_{2} } \right|_{\tilde{x}_{l} ^{*} } \left(x_{2} -\tilde{x}_{2}^{*} \right)+...+\left. \frac{\mu _{k} \left(f_{k} \left(\tilde{x}_{l}^{*} \right)\right)}{\partial x_{n} } \right|_{\tilde{x}_{l} ^{*} } \left(x_{n} -\tilde{x}_{n}^{*} \right)\right] 
\end{equation} 
\end{enumerate}

Here, functions $\tilde{\mu }_{k} \left(f_{k} \left(x\right)\right)$approximate the nonlinear functions $\mu _{k} \left(f_{k} \left(x\right)\right)$ around the maximal solution $\tilde{x}_{l}^{*} =\left(\tilde{x}_{1}^{*} ,\tilde{x}_{2}^{*} ,...,\tilde{x}_{n}^{*} \right).$  So, Taylor series approach generally provides a relatively good approximation to a differentiable function but only around a given point, and not over the entire domain.

\subsection{ A Fuzzy Goal Programming Model to Multiobjective Linear Programming Problem}

The FGP approach was originally introduced by Zimmermann [4] in 1978. He employed the concept of membership functions. Tiwari et al. [35] suggested a weighted additive model that associates each goal's weight into the objective function, where weights show the relative importance of the fuzzy goals. Afterward, Mohamed [36] suggested a kind of fuzzy goal, which is introduced in the general form of FGP model. In [43], Mohamed's approach used to present a FGP approach for solving multiobjective programming problems and then Gupta and Bhattacharjee [40] formulated two FGP model for solving multiobjective programming problems. 

According to paper [36], the highest degree of membership function is 1 and therefore, the nonlinear membership functions in \eqref{GrindEQ__15_} and \eqref{GrindEQ__17_} can be constructed as the following nonlinear membership goals;
\begin{equation} \label{GrindEQ__22_} 
\mu _{k} \left(f_{k} \left(x\right)\right)\cong \frac{f_{k} \left(x\right)-L_{k} }{\left(U_{k} -L_{k} \right)} +d_{k}^{-} -d_{k}^{+} =1,k=1,2,..,l. 
\end{equation} 
\begin{equation} \label{GrindEQ__23_} 
\mu _{k} \left(f_{k} \left(x\right)\right)=\frac{U_{k} -f_{k} \left(x\right)}{U_{k} -L_{k} } +d_{k}^{-} -d_{k}^{+} =1,k=1,2,..,l 
\end{equation} 
where $d_{k}^{-} \left(\ge 0\right)$and $d_{k}^{+} \left(\ge 0\right)$ represent the negative and positive deviations from the aspired levels, respectively. In addition, any positive deviation from 1 shows the full attainment of the membership value. Hence to reach the aspired levels of the fuzzy goal, it is sufficient to minimize its negative deviational variable from 1. At the same time, presentation of both deviation variables in the membership goal is unneeded, and the positive deviational variables are not necessary [40]. Thus the above membership goals can be written as follows:
\begin{equation} \label{GrindEQ__24_} 
\mu _{k} \left(f_{k} \left(x\right)\right)+d_{k}^{-} =1, k=1,2,..,l 
\end{equation} 
Here, the membership goals as defined in \eqref{GrindEQ__24_} are naturally nonlinear when the objective functions are nonlinear, and this may generate computational difficulties in the solution procedure of nonlinear problems. Therefore, by using Taylor series approach, the nonlinear membership goal in \eqref{GrindEQ__24_} can be written as the following linear function:
\begin{equation} \label{GrindEQ__25_} 
\tilde{\mu }_{k} \left(f_{k} \left(x\right)\right)+d_{k}^{-} -d_{k}^{+} =1, k=1,2,..,l 
\end{equation} 
Now let us consider the min-max form of fuzzy programming \eqref{GrindEQ__20_}. If we put $\beta =\, \left(1-\lambda \right),$then we obtain the following equivalent fuzzy linear programming model:
\begin{equation} \label{GrindEQ__26_} 
\begin{array}{l} {Min\, \, \, \beta } \\ {\left\{\begin{array}{l} {\, \beta \ge 1-\tilde{\mu }_{k} \left(f_{k} \left(x\right)\right),\, \, k=1,2,..,l,} \\ {\underline{x}_{l}^{} \le x_{l}^{} \le \bar{x}_{l}^{} ,\, \, \, l=1,2,...,n,} \\ {\beta \le 1,\, \, \beta \ge 0\, } \end{array}\right. } \end{array} 
\end{equation} 
where $\underline{x}_{l}^{} \le x_{l}^{} \le \bar{x}_{l}^{} $ denotes that the limits of decision variables derived from the individual optimal solutions of each objective.

\noindent To formulate the above fuzzy problem as a FGP model, the negative deviational variables in \eqref{GrindEQ__26_} can be defined as:
\[d_{k}^{-} =Max\left\{0,1-\tilde{\mu }_{k} \left(f_{k} \left(x\right)\right)\right\}\] 
Thus, we obtain $\beta \ge d_{k}^{-} ,$ where $\tilde{\mu }_{k} \left(f_{k} \left(x\right)\right)+d_{k}^{-} =1,$$k=1,2,..,l.$ 

\noindent Then, an equivalent linear FGP model for problem \eqref{GrindEQ__13_} can be developed as follows:
\begin{equation} \label{GrindEQ__27_} 
\begin{array}{l} {Min\, \, \beta } \\ {\left\{\begin{array}{l} {\tilde{\mu }_{k} \left(f_{k} \left(x\right)\right)+d_{k}^{-} =1,\, \, \, k=1,2,..,l,} \\ {\beta \ge d_{k}^{-} ,\, k=1,2,..,l,\, \, } \\ {\underline{x}_{l}^{} \le x_{l}^{} \le \bar{x}_{l}^{} ,\, \, \, l=1,2,...,n.} \\ {d_{k}^{-} \ge 0,\, \, \, k=1,2,..,l,} \\ {\beta \ge 0,\, \, \beta \le 1} \end{array}\right. } \end{array} 
\end{equation} 
where $d_{k}^{-} \left(\ge 0\right)$represent the negative deviations from the aspired levels, respectively.

\subsection{ Interactive Fuzzy Goal Programming Approach Based on Taylor Series for MNLOPP with IT2 FNs}

\noindent In this section, an interactive fuzzy goal programming algorithm is presented to achieve the highest degree for the membership functions. 

\noindent The complete suggested solution procedures can be summarized as follows.

\noindent \textbf{Step 1} Construct the mathematical model of MONLPP with IT2 FNs \eqref{GrindEQ__12_}.

\noindent \textbf{Step 2} By using the expected value function as defined in \eqref{GrindEQ__9_}, obtain the corresponding crisp MONPP.

\noindent \textbf{Step 3} Solve the MONPP as a single objective problem, considering each time only one objective as the objective function and ignoring all others. 

\noindent \textbf{Step 4 }Compute the value of each objective function at each solution derived in Step 3. Then define the feasible region by the lower and upper limits of decision variables.

\noindent \textbf{Step 5} From Step 4, determine the upper and lower tolerance limits of each objective function.

\noindent \textbf{Step 6} Construct the membership functions as defined in \eqref{GrindEQ__15_} and \eqref{GrindEQ__17_} for each objective.   

\noindent \textbf{Step 7} Maximize each nonlinear membership functions under the feasible region derived in Step 4, individually and then determine the maximal solutions $\tilde{x}_{l}^{*} =\left(\tilde{x}_{1}^{*} ,\tilde{x}_{2}^{*} ,...,\tilde{x}_{n}^{*} \right)$for each nonlinear membership function. 

\noindent \textbf{Step 8} Linearize each nonlinear membership function using Taylor series at the maximal solutions $\tilde{x}_{l}^{*} =\left(\tilde{x}_{1}^{*} ,\tilde{x}_{2}^{*} ,...,\tilde{x}_{n}^{*} \right).$

\noindent \textbf{Step 9} Construct the FGP model as formulated in \eqref{GrindEQ__27_}, then solve it to obtain the optimal solution.

\noindent \textbf{Step 10} If the decision-maker is satisfied by the current solution, in Step 9, go to Step 11, else go to Step 12.

\noindent \textbf{Step 11} The current solution is the optimal solution for MONPP with IT2 FNs.

\noindent \textbf{Step 12} Compare the lower (upper) tolerance limit of each objective with the new value of the objective function. If the new value is higher (lower) than the lower tolerance limit, take this as a new lower (upper) tolerance limit. If else, hold the old one as is and then go to step \eqref{GrindEQ__5_}.

\section{ Numerical Examples}

\textbf{Example 1}

A manufacturing factory is going to produce 3 kinds of products A; B and C in a period (say one month). The production of A; B and C require three kinds of resources $R_{1} ,$$R_{2} $ and $R_{3} $. Here, to deal with uncertainties possessing doubt, let us consider that all parameters of the problem are IT2 FNs. Thus, the requirements of each kind of resource to produce each product A are $a_{11} =\left(\left({\rm 3,}\, {\rm 3,}\, {\rm 4,}\, {\rm 5;}\, {\rm 0.90,}\, {\rm 0.91}\right),\right. $$\left. \left({\rm 4,}\, {\rm 4,}\, {\rm 5,}\, {\rm 6;\; 0.}\, {\rm 92,0.}\, {\rm 93}\right)\right);$ $a_{12} =\left(\left({\rm 3,}\, {\rm 5,}\, {\rm 5,}\, {\rm 7;}\, {\rm 0.90,}\, {\rm 0.98}\right)\right. ,$$\left. \left({\rm 2,}\, {\rm 4,}\, {\rm 4,}\, {\rm 5;}\, {\rm 0.92,}\, {\rm 0.97}\right)\right);$$a_{13} =\left(\left({\rm 2,}\, {\rm 4,}\, {\rm 4,}\, {\rm 5;}\, {\rm 0.90,}\, {\rm 0.91}\right){\rm ,}\right. $$\left. \left({\rm 3,}\, {\rm 4,}\, {\rm 5,}\, {\rm 5;}\, {\rm 0.92,}\, {\rm 0.93}\right)\right);$ units, respectively. To produce each product B, the respective requirements are $a_{21} =\left(\left({\rm 3,5,5,7;0.90,0.98}\right)\right. ,$$\left. \left({\rm 2,4,4,5;0.92,0.97}\right)\right);$$a_{22} =\left(\left({\rm 3,3,4,5;0.90,0.91}\right){\rm ,}\right. $$\left. \left({\rm 4,4,5,6,0.92,0.93}\right)\right);$ $a_{23} =\left(\left({\rm 3,3,4,5;0.90,0.91}\right){\rm ,}\right. $$\left. \left({\rm 4,4,5,6;0.92,0.93}\right)\right);$ units  and that for each unit of C are around $a_{31} =\left(\left({\rm 2,4,4,5;0.90,0.91}\right){\rm ,}\right. $$\left. \left({\rm 3,4,5,5;0.92,0.93}\right)\right);$ $a_{32} =\left(\left({\rm 3,3,4,5;0.90,0.91}\right){\rm ,}\right. $$\left. \left({\rm 4,4,5,6;0.92,0.93}\right)\right);$$a_{33} =\left(\left({\rm 2,4,4,5;0.90,0.91}\right){\rm ,}\right. $$\left. \left({\rm 3,4,5,5;0.92,0.93}\right)\right);$ units. The planned existing resource of $R_{1} ,$$R_{2} $ are around $b_{1} =\left(\left({\rm 80,95,70,90;0.96,0.99}\right){\rm ,}\right. $$\left. \left({\rm 90,80,100,110;0.97,0.99}\right)\right);$$b_{2} =\left(\left({\rm 90,50,70,70;0.95,0.98}\right){\rm ,}\right. $$\left. \left({\rm 90,80,80,90;0.97,0.99}\right)\right)$ units respectively. But there is additional safety store of materials, which are administrated by the manager. For better quality of the products, at least, $b_{3} =\left(\left({\rm 50,60,60,70;0.95,0.99}\right){\rm ,\; }\left({\rm 50,60,60,70;0.94,0.99}\right)\right)$  units of resource $R_{3} $ has to be employed. In addition, the conjectural time requirements in producing each unit of products are $\tilde{t}_{1} ,$$\tilde{t}_{2} ,$ and $\tilde{t}_{3} $ h respectively. 

\noindent Assuming that the planned production quantities of A; B and C are $x_{1} ;x_{2} ;x_{3}^{} $respectively. Moreover, assuming that unit cost and sale's price of product A, B and C are $UC_{1} =\tilde{c}_{1} ,$ $UC_{2} =\tilde{c}_{2} ,$$UC_{3} =\tilde{c}_{3} $and $US_{1} =\frac{\tilde{s}_{1} }{x_{1}^{{\raise0.7ex\hbox{$ 1 $}\!\mathord{\left/ {\vphantom {1 a_{1} }} \right. \kern-\nulldelimiterspace}\!\lower0.7ex\hbox{$ a_{1}  $}} } } ,$ $US_{2} =\frac{\tilde{s}_{2} }{x_{2}^{{\raise0.7ex\hbox{$ 1 $}\!\mathord{\left/ {\vphantom {1 a_{2} }} \right. \kern-\nulldelimiterspace}\!\lower0.7ex\hbox{$ a_{2}  $}} } } ,$ $US_{3} =\frac{\tilde{s}_{3} }{x_{1}^{{\raise0.7ex\hbox{$ 1 $}\!\mathord{\left/ {\vphantom {1 a_{3} }} \right. \kern-\nulldelimiterspace}\!\lower0.7ex\hbox{$ a_{3}  $}} } } $respectively, where $a_{1} =2;a_{2} =2;a_{3}^{} =3$are real numbers. Here, the decision maker expects to maximize whole profit and minimize integral time requirement. 

\noindent (Step 1): This problem can be formulated as follows: 
\begin{equation} \label{GrindEQ__28_} 
\begin{array}{l} {Max\, f_{2} \left(x\right)=\tilde{s}_{1} x_{1}^{1-{\raise0.7ex\hbox{$ 1 $}\!\mathord{\left/ {\vphantom {1 a_{1} }} \right. \kern-\nulldelimiterspace}\!\lower0.7ex\hbox{$ a_{1}  $}} } -\tilde{c}_{1} x_{1} +\tilde{s}_{2} x_{2}^{1-{\raise0.7ex\hbox{$ 1 $}\!\mathord{\left/ {\vphantom {1 a_{2} }} \right. \kern-\nulldelimiterspace}\!\lower0.7ex\hbox{$ a_{2}  $}} } -\tilde{c}_{2} x_{2} +\tilde{s}_{3} x_{3}^{1-{\raise0.7ex\hbox{$ 1 $}\!\mathord{\left/ {\vphantom {1 a_{3} }} \right. \kern-\nulldelimiterspace}\!\lower0.7ex\hbox{$ a_{3}  $}} } -\tilde{c}_{3} x_{3} } \\ {Min\, f_{2} \left(x\right)=\tilde{t}_{1} x_{1} +\tilde{t}_{2} x_{2} +\tilde{t}_{3} x_{3} ,} \\ {s.t.\left\{\begin{array}{l} {\tilde{a}_{11} x_{1} x_{2}^{} +\tilde{a}_{12} x_{2} +\tilde{a}_{13} x_{3}^{2} \le \tilde{b}_{1} ,} \\ {\tilde{a}_{21} x_{1}^{2} +\tilde{a}_{22} x_{1}^{} x_{2} +\tilde{a}_{23} x_{3} \le \tilde{b}_{2} ,} \\ {\tilde{a}_{31} x_{1} +\tilde{a}_{32} x_{2} +\tilde{a}_{33} x_{3} \ge \tilde{b}_{3} } \\ {x_{1} ,x_{2} ,x_{3} \ge 0} \end{array}\right. } \end{array} 
\end{equation} 
where $\tilde{s}_{1} =\left(\left({\rm 20,22,24,27;0.95,0.98}\right),\right. $$\left. \left({\rm 21,23,25,26;0.97,0.99}\right)\right),$ $\tilde{s}_{2} =\left(\left({\rm 21,23,24,28;0.94,0.99}\right){\rm ,}\right. $$\left. \left({\rm 22,23,25,26;0.95,0.97}\right)\right),$$\tilde{s}_{3} =\left(\left({\rm 22,23,24,26;0.94,0.97}\right)\right. ,$$\left. \left({\rm 22,24,25,26;0.95,0.97}\right)\right);$$\tilde{c}_{1} =\left(\left({\rm 1,3,3,4;0.90,0.91}\right){\rm ,}\right. $$\left. \left({\rm 1,2,4,5;0.92,0.93}\right)\right),$$\tilde{c}_{2} =\left(\left({\rm 2,3,5,5;0.91,0.94}\right){\rm ,}\right. $$\left(\left({\rm 2,3,6,8;0.93,0.95}\right)\right),$$\tilde{c}_{3} =\left(\left({\rm 2,4,4,5;0.90,0.91}\right){\rm ,}\right. $$\left. \left({\rm 3,4,5,5;0.92,0.93}\right)\right),$$\tilde{t}_{1} =\left(\left({\rm 2,3,4,5;0.95,0.99}\right)\right. ,$$\left. \left({\rm 1,2,3,3;0.92,0.97}\right)\right),$$\tilde{t}_{2} =\left(\left({\rm 3,4,5,6;0.96,0.98}\right),\right. $$\left. \left({\rm 1,2,3,3;0.95,0.96}\right)\right),$  $\tilde{t}_{3} =\left(\left({\rm 3,3,4,5;0.90,0.91}\right),\right. $$\left. \left({\rm 4,4,5,6;0.92,0.93}\right)\right)$ are the estimated coefficient parameters.

\noindent (Step 2): Employing the expected value function in \eqref{GrindEQ__9_}, problem \eqref{GrindEQ__28_} is transformed into an equivalent crisp MONLPP as follows:
\begin{equation} \label{GrindEQ__29_} 
\begin{array}{rcl} {Max\, \, f_{1} \left(x\right)} & {=} & {22.854x_{1} ^{(1/2)} {\rm -2.631}x_{1} +{\rm 23.100}x_{2} ^{(1/2)} -{\rm 3.963}x_{2} +{\rm 22.980}x_{3} ^{(2/3)} {\rm -3.660}x_{3} } \\ {} \end{array} 
\end{equation} 
(Step 3): Then, problem \eqref{GrindEQ__29_} is solved as a single objective problem and the individual maximum and minimum solutions for each objective are given in Table 1.

\begin{center}{Table 1: The individual optimal solutions for each objective}

\begin{tabular}{|p{0.7in}|p{0.7in}|p{0.7in}|p{0.7in}|p{0.7in}|} \hline 
$x$ & $Max\, \, f_{1} \left(x\right)$ & $Min\, \, f_{1} \left(x\right)$ & $Max\, \, f_{2} \left(x\right)$ & $Min\, \, f_{2} \left(x\right)$ \\ \hline 
$x_{1}^{} $ & $0.203$ & $0$ & $0$ & $0.457$ \\ \hline 
$x_{2}^{} $ & $12.344$ & $15.237$ & $20.862$ & $14.808$ \\ \hline 
$x_{3}^{} $ & $2.703$ & $0$ & $0.607$ & $0$ \\ \hline 
\end{tabular}
\end{center}

\noindent (Step 4): From Table 1, the objective function values are determined as follows:

\noindent $29.785\le f_{1} \le 76.694$and $54.699\le f_{2} \le 77.653$,

\noindent In addition, the nonlinear constrained region is reduced to the following inequalities:
\[0\le x_{1} \le 0.457, 12.344\le x_{2} \le 20.862, 0\le x_{3} \le 2.703\] 
 (Step 5): Then, aspiration level of each objective function is $s_{1} =76.694$ and $s_{2} =54.699,$respectively; The lower tolerance limit of the first objective function is $29.785;$The upper tolerance limit of  the second objective is $77.653.$

\noindent (Step 6): Then the membership functions can be constructed as follows based on \eqref{GrindEQ__15_} and \eqref{GrindEQ__17_}:
\begin{equation} \label{GrindEQ__30_} 
\mu _{1} \left(f_{1} \left(x\right)\right)\cong \left\{\begin{array}{l} {1\, \, \, \, \, \, \, \, \, \, \, \, \, \, \, \, \, \, \, \, \, \, \, \, \, \, \, \, \, \, \, \, \, \, \, \, \, \, \, \, f_{1} \left(x\right)\ge 77.653,} \\ {\frac{f_{1} \left(x\right)-29.785}{\left(77.653-29.785\right)} ,\, \, 29.785\le f_{1} \left(x\right)\le 77.653,\, \, \, \, \, \, \, } \\ {0\, \, \, \, \, \, \, \, \, \, \, \, \, \, \, \, \, \, \, \, \, \, \, \, \, \, \, \, \, \, \, \, \, \, \, \, \, \, \, 29.785\ge f_{1} \left(x\right)} \end{array}\right.  
\end{equation} 
\[=0.487\sqrt{x_{1} } -0.056x_{1} +0.492\sqrt{x_{2} } -0.085x_{2} +0.490x_{3}^{\left(2/3\right)} -0.078x_{3} -0.635\] 
\begin{equation} \label{GrindEQ__31_} 
\mu _{2} \left(f_{2} \left(x\right)\right)=\left\{\begin{array}{l} {1\, \, \, \, \, \, \, \, \, \, \, \, \, \, \, \, \, \, \, \, \, \, \, \, \, \, \, \, \, \, \, \, \, \, \, \, \, \, \, \, \, \, \, \, \, \, \, \, f_{2} \left(x\right)\le 54.699,} \\ {\frac{77.653-f_{2} \left(x\right)}{77.653-54.699} ,\, \, \, \, 54.699\le f_{2} \left(x\right)\le 77.653} \\ {0\, \, \, \, \, \, \, \, \, \, \, \, \, \, \, \, \, \, \, \, \, \, \, \, \, \, \, \, \, \, \, \, \, \, \, \, \, \, \, \, \, \, \, \, \, \, \, 77.653\le f_{2} \left(x\right)} \end{array}\right.  
\end{equation} 
\[=-0.120x_{1} -0.157x_{2} -0.169x_{3} +3.383\] 
(Step 7): The maximal solution of membership function \eqref{GrindEQ__30_} under the constraints is determined as follows:
\[\tilde{x}_{}^{*} =\left(\tilde{x}_{1}^{*} =0.458,\, \, \tilde{x}_{2}^{*} ={\rm 12.344,}\, \tilde{x}_{3}^{*} =2.703\right).\] 
 (Step 8): Thereby the nonlinear membership function in \eqref{GrindEQ__30_} is converted into the linear functions using Taylor series approach around the maximal solution $\tilde{x}_{}^{*} =\left(\tilde{x}_{1}^{*} =0.458,\, \, \tilde{x}_{2}^{*} ={\rm 12.344,}\, \tilde{x}_{3}^{*} =2.703\right).$ 

\noindent Thus, the linear membership function is as follows:
\[\begin{array}{rcl} {} & {} & {\tilde{\mu }_{1} \left(f_{1} \left(x\right)\right)\cong \left[\left. \frac{\mu _{1} \left(f_{1} \left(\tilde{x}_{}^{*} \right)\right)}{\partial x_{1} } \right|_{\tilde{x}^{*} } \left(x_{1} -0.457\right)+\left. \frac{\mu _{1} \left(f_{1} \left(\tilde{x}_{}^{*} \right)\right)}{\partial x_{2} } \right|_{\tilde{x}^{*} } \left(x_{2} -12.344\right)+\left. \frac{\mu _{3} \left(f_{3} \left(\tilde{x}_{}^{*} \right)\right)}{\partial x_{3} } \right|_{\tilde{x}^{*} } \left(x_{3} -2.703\right)\right]} \\ {} & {=} & {0{\rm .304}x_{1} -0{\rm .014}x_{2} +0.156x_{3} +{\rm 0.712}} \end{array}\] 
(Step 9): Consequently, the proposed linear FGP model is constructed as follows based on model \eqref{GrindEQ__26_};
\begin{equation} \label{GrindEQ__32_} 
\begin{array}{l} {Min\, \, \beta } \\ {s.t.\left\{\begin{array}{l} {0{\rm .304}x_{1} -0{\rm .014}x_{2} +0.156x_{3} +{\rm 0.712}+d_{1}^{-} =1,\, } \\ {-0.120x_{1} -0{\rm .157}x_{2} +0.169x_{3} +3{\rm .383}+d_{2}^{-} =1,} \\ {0\le x_{1} \le 0.457,} \\ {12.344\le x_{2} \le 20.862,} \\ {0\le x_{3} \le 2.703,} \\ {\beta \ge d_{1}^{-} ,\beta \ge d_{2}^{-} ,} \\ {\beta \ge 0,\, \, \beta \le 1,} \\ {d_{1}^{-} ,\, \, d_{2}^{-} \ge 0.} \end{array}\right. } \end{array} 
\end{equation} 
The above problem is solved by using Maple 18.02 optimization toolbox and then the optimal solution of problem \eqref{GrindEQ__32_} are as follows: $\beta =0,\, \, d_{1}^{-} =0,\, \, d_{2}^{-} =0,\, $$x_{1} =0.458,\, \, \, x_{2} ={\rm 12.710},\, \, \, x_{3} =1.946,\, $with the objective function values as:$f_{1} ={\rm 74.938},\, \, \, f_{2} ={\rm 54.699}.$Also the membership values are as follows:$\mu _{1} \left(f_{1} \right)=0.963,$$\, \mu _{2} \left(f_{2} \right)=1.00.$ 

\noindent (Step 10): Let the decision maker be satisfied by the optimal solution $\left(x_{1} =0.458,\, \, \, x_{2} ={\rm 12.710},\, \, \, x_{3} =1.946\right),\, $and then the proposed algorithm is stopped at step 11.

\noindent In order to demonstrate the performance of the suggested procedure, the above numerical example is solved by using different fuzzy goal programming approaches. In these approaches, the numerical weights are considered as follows:
\[w_{k}^{} =\frac{1}{\varepsilon } ,k=1,2,..,l\] 
where $\varepsilon $ represent the tolerant intervals. Comparative results are given in Table 3.

\begin{center}{Tablo 2 Comparison of results by different approaches}

\begin{tabular}{|p{0.4in}|p{1.2in}|p{1.2in}|p{1.2in}|p{1.2in}|} \hline 
$ $  & The proposed approach \newline$\left(0.458,12.710,1.946\right)$  & Mohamed's approach [36]\newline $\left(0.348, 13.677, 1.549\right)$ & Model I in [40]\newline  $\left(0.341,13.782, 1.418\right)$ & Model II in [40]\newline $\left(0.315, 12.607, 2.334\right)$ \\ \hline 
$f_{1} $ & $74.938$ & $68.894$ & $67.404$ & $75.950$ \\ \hline 
$\, f_{2} $ & $54.699$ & $56.342$ & $56.191$ & $55.443$ \\ \hline 
$\mu _{1} \left(f_{1} \right)$ & $0.963$ & $0.834$ & $0.801$ & $0.984$ \\ \hline 
$\mu _{2} \left(f_{2} \right)$ & $1.00$ & $0.928$ & $0.935$ & $0.968$ \\ \hline 
\end{tabular}
\end{center}
From Table 2, all of the sums of the membership values generated by the suggested procedure is greater than that generated by the approaches in [36] and [40].  

\textbf{Example 2}

In order to further verify the correctness of the suggested procedure, let us consider the following data Table 
\begin{center}{Tablo 3 Estimated trapezoidal IT2 FNs in the problem}
\newline
\begin{tabular}{|p{3.5in}|p{3.5in}|} \hline 
\textit{Fuzzy unit cost in the first objective\newline }$\tilde{s}_{1} =\left(\left({\rm 80,95,70,90;0.96,0.99}\right){\rm ,}\right. $$\left. \left({\rm 90,80,100,110;0.97,0.99}\right)\right)$\newline $\tilde{s}_{2} =\left(\left({\rm 50,60,60,70;0.95,0.99}\right){\rm ,\; }\left({\rm 50,60,60,70;0.94,0.99}\right)\right)$\newline $\tilde{s}_{3} =\left(\left({\rm 90,50,70,70;0.95,0.98}\right){\rm ,}\right. $$\left. \left({\rm 90,80,80,90;0.97,0.99}\right)\right)$\newline $\tilde{c}_{1} =\left(\left({\rm 3,5,5,7;0.90,0.98}\right)\right. ,$$\left. \left({\rm 2,4,4,5;0.92,0.97}\right)\right)$\newline $\tilde{c}_{2} =\left(\left({\rm 2,3,5,5;0.91,0.94}\right){\rm ,}\right. $$\left(\left({\rm 2,3,6,8;0.93,0.95}\right)\right)$\newline $\tilde{c}_{3} =\left(\left({\rm 3,}\, {\rm 3,}\, {\rm 4,}\, {\rm 5;}\, {\rm 0.90,}\, {\rm 0.91}\right),\right. $$\left. \left({\rm 4,}\, {\rm 4,}\, {\rm 5,}\, {\rm 6;\; 0.}\, {\rm 92,0.}\, {\rm 93}\right)\right)$\newline \textit{Fuzzy unit time in the second objective\newline }$\tilde{t}_{1} =\left(\left({\rm 3,5,5,7;0.90,0.98}\right)\right. ,$$\left. \left({\rm 2,4,4,5;0.92,0.97}\right)\right)$\newline $\tilde{t}_{2} =\left(\left({\rm 3,4,5,6;0.96,0.98}\right),\right. $$\left. \left({\rm 1,2,3,3;0.95,0.96}\right)\right)$\newline $\tilde{t}_{3} =\left(\left({\rm 2,3,4,5;0.95,0.99}\right)\right. ,$$\left. \left({\rm 1,2,3,3;0.92,0.97}\right)\right)$\newline \textit{Fuzzy requirements of resources\newline }$\tilde{a}_{11} =\left(\left({\rm 2,3,4,5;0.95,0.99}\right)\right. ,$$\left. \left({\rm 1,2,3,3;0.92,0.97}\right)\right)$\newline $\tilde{a}_{12} =\left(\left({\rm 2,3,4,5;0.95,0.99}\right)\right. ,$$\left. \left({\rm 1,2,3,3;0.92,0.97}\right)\right)$\newline  & $\tilde{a}_{13} =\left(\left({\rm 3,4,5,6;0.96,0.98}\right),\right. $$\left. \left({\rm 1,2,3,3;0.95,0.96}\right)\right)$\newline $\tilde{a}_{21} =\left(\left({\rm 2,3,4,5;0.95,0.99}\right)\right. ,$$\left. \left({\rm 1,2,3,3;0.92,0.97}\right)\right)$\newline $\tilde{a}_{22} =\left(\left({\rm 3,}\, {\rm 3,}\, {\rm 4,}\, {\rm 5;}\, {\rm 0.90,}\, {\rm 0.91}\right),\right. $$\left. \left({\rm 4,}\, {\rm 4,}\, {\rm 5,}\, {\rm 6;\; 0.}\, {\rm 92,0.}\, {\rm 93}\right)\right)$\newline $\tilde{a}_{23} =\left(\left({\rm 3,}\, {\rm 3,}\, {\rm 4,}\, {\rm 5;}\, {\rm 0.90,}\, {\rm 0.91}\right),\right. $$\left. \left({\rm 4,}\, {\rm 4,}\, {\rm 5,}\, {\rm 6;\; 0.}\, {\rm 92,0.}\, {\rm 93}\right)\right)$\newline $\tilde{a}_{31} =\left(\left({\rm 3,4,5,6;0.96,0.98}\right),\right. $$\left. \left({\rm 1,2,3,3;0.95,0.96}\right)\right)$\newline $\tilde{a}_{32} =\left(\left({\rm 3,}\, {\rm 3,}\, {\rm 4,}\, {\rm 5;}\, {\rm 0.90,}\, {\rm 0.91}\right),\right. $$\left. \left({\rm 4,}\, {\rm 4,}\, {\rm 5,}\, {\rm 6;\; 0.}\, {\rm 92,0.}\, {\rm 93}\right)\right)$\newline $\tilde{a}_{33} =\left(\left({\rm 3,}\, {\rm 3,}\, {\rm 4,}\, {\rm 5;}\, {\rm 0.90,}\, {\rm 0.91}\right),\right. $$\left. \left({\rm 4,}\, {\rm 4,}\, {\rm 5,}\, {\rm 6;\; 0.}\, {\rm 92,0.}\, {\rm 93}\right)\right)$\newline \textit{Available fuzzy resources\newline }$\tilde{b}_{1} =\left(\left({\rm 20,22,24,27;0.95,0.98}\right),\right. $$\left. \left({\rm 21,23,25,26;0.97,0.99}\right)\right)$ $\tilde{b}_{2} =\left(\left({\rm 22,23,24,26;0.94,0.97}\right)\right. ,$$\left. \left({\rm 22,24,25,26;0.95,0.97}\right)\right)$\newline $\tilde{b}_{3} =\left(\left({\rm 21,23,24,28;0.94,0.99}\right){\rm ,}\right. $$\left. \left({\rm 22,23,25,26;0.95,0.97}\right)\right)$\newline and\newline $a_{1} =2;a_{2} =2;a_{3}^{} =3$\textit{} \\ \hline 
\end{tabular}
\end{center}
Using the expected value function as defined in \eqref{GrindEQ__9_} (Step 2), all the IT2 FNs in data Table 3 can be converted into the crisp numbers as in the following data Table 4.

\begin{center}{Table 4 The expected values of estimated trapezoidal IT2 FNs for the problem}

\begin{tabular}{|p{2.5in}|p{2.5in}|} \hline 
\textit{Unit cost in the first objective\newline }$s_{11} =87.364$$s_{12} =59.259$$s_{13} ={\rm 75.369}$\newline $c_{11} =4.123$$\tilde{c}_{12} ={\rm 3.963}$$c_{13} =3.889$\newline \textit{Unit time in the second objective\newline }$t_{11} ={\rm 4.123}$$\tilde{t}_{12} =3.609$$\tilde{t}_{13} =2.753,$ & \textit{Requirements of each of resources\newline }$a_{11} ={\rm 2.753},$$a_{12} ={\rm 2.753},$$a_{13} ={\rm 3.609}$\newline $a_{21} ={\rm 2.753}$$a_{22} ={\rm 3.889}$$a_{23} ={\rm 3.889}$\newline $a_{31} ={\rm 3.609}$$a_{32} ={\rm 3.889},$$a_{33} ={\rm 3.889}$\textit{\newline Available resources\newline }$b_{1} ={\rm 22.854},$$b_{2} ={\rm 22.980}$$b_{3} ={\rm 23.1}$\textit{} \\ \hline 
\end{tabular}
\end{center}

From table 4, the considered MONLPP with IT2 FNs is transformed as follows based on model \eqref{GrindEQ__13_}:
\begin{equation} \label{GrindEQ__33_} 
\begin{array}{rcl} {Max\, \, f_{1} \left(x\right)} & {=} & {87.364\sqrt{x_{1} } {\rm -4.123}x_{1} +{\rm 59.259}\sqrt{x_{2} } -{\rm 3.963}x_{2} +{\rm 75.369}\sqrt[{3}]{x_{2}^{2} } {\rm -3.889}x_{3} } \\ {} \end{array} 
\end{equation} 
(Step 3)The individual maximum and minimum solutions are summarized in Table 5.

\begin{center}{Table 5 The individual minimum and maximum solutions }
\begin{tabular}{|p{0.2in}|p{0.7in}|p{0.7in}|p{0.7in}|p{0.7in}|} \hline 
$x$ & $Max\, \, f_{1} \left(x\right)$ & $Min\, \, f_{1} \left(x\right)$ & $Max\, \, f_{2} \left(x\right)$ & $Min\, \, f_{2} \left(x\right)$ \\ \hline 
$x_{1}^{} $ & $1.051$ & $0$ & $0$ & $1.035$ \\ \hline 
$x_{2}^{} $ & $3.209$ & $5.940$ & $8.053$ & $4.980$ \\ \hline 
$x_{3}^{} $ & $1.756$ & $0$ & $0.436$ & $0$ \\ \hline 
\end{tabular}
\end{center}
\noindent (Step 4) The objective function values are as follows:
\[120.885\le f_{1} \le 281.523;20.820\le f_{2} \le 30.858.\] 
Then the feasible region with the bounded decision variables are as follows:
\[0\le x_{1} \le 1.051, 3.209\le x_{2} \le 8.053, 0\le x_{3} \le 1.756\] 
(Step 5) Aspiration level for each objective function is $s_{1} =281.523$and $s_{2} =20.820$respectively; The lower and upper tolerance limits of the objective functions are $120.885$and $30.858,$ respectively. 

\noindent (Step 6)Then the membership functions are formulated as:
\[\mu _{1} \left(f_{1} \left(x\right)\right)\cong \left\{\begin{array}{l} {1\, \, \, \, \, \, \, \, \, \, \, \, \, \, \, \, \, \, \, \, \, \, \, \, \, \, \, \, \, \, \, \, \, \, \, \, \, \, \, \, \, \, \, \, f_{1} \left(x\right)\ge 281.523,} \\ {\frac{f_{1} \left(x\right)-120.885}{\left(281.523-120.885\right)} ,\, \, 120.885\le f_{1} \left(x\right)\le 281.523,\, \, \, \, \, \, \, } \\ {0\, \, \, \, \, \, \, \, \, \, \, \, \, \, \, \, \, \, \, \, \, \, \, \, \, \, \, \, \, \, \, \, \, \, \, \, \, \, \, \, \, \, \, \, 120.885\ge f_{1} \left(x\right)} \end{array}\right. \] 
\[=0.544\sqrt{x_{1} } -0.026x_{1} +0.369\sqrt{x_{2} } -0.025x_{2} +0.469\sqrt[{3}]{x_{3}^{2} } -0.024x_{3} -0.753\] 
\[\mu _{3} \left(f_{3} \left(x\right)\right)=\left\{\begin{array}{l} {1\, \, \, \, \, \, \, \, \, \, \, \, \, \, \, \, \, \, \, \, \, \, \, \, \, \, \, \, \, \, \, \, \, \, \, \, \, \, \, \, \, \, \, \, \, \, \, \, f_{3} \left(x\right)\le 20.820,} \\ {\frac{30.858-f_{3} \left(x\right)}{30.858-20.820} ,\, \, \, \, \, \, \, \, \, \, \, \, \, \, 20.820\le f_{3} \left(x\right)\le 30.858} \\ {0\, \, \, \, \, \, \, \, \, \, \, \, \, \, \, \, \, \, \, \, \, \, \, \, \, \, \, \, \, \, \, \, \, \, \, \, \, \, \, \, \, \, \, \, \, \, \, 30.858\le f_{3} \left(x\right)} \end{array}\right. \] 
\[=-0.274x_{1} -0.360x_{2} -0.411x_{3} +3.074\] 
The maximal solutions for the nonlinear membership function $\mu _{1} \left(f_{1} \left(x\right)\right)$ under the constraints (Step 7) is: 
\[\tilde{x}_{}^{*} =\left(\tilde{x}_{1}^{*} =1.051,\, \, \tilde{x}_{2}^{*} ={\rm 8.053,}\, \tilde{x}_{3}^{*} =1.756\right).\] 
(Step 8) Then, the nonlinear membership function $\mu _{1} \left(f_{1} \left(x\right)\right)$ is converted into an equivalent linear function using Taylor series approach around its maximal solution $\tilde{x}_{}^{*} =\left(\tilde{x}_{1}^{*} =1.051,\, \, \tilde{x}_{2}^{*} ={\rm 8.053,}\, \tilde{x}_{3}^{*} =1.756\right).$ Thereby, the linear function is determined as:
\[\tilde{\mu }_{1} \left(f_{1} \left(x\right)\right)=0{\rm .240}x_{1} -0{\rm .040}x_{2} +0.235x_{3} +{\rm 0.277.}\] 
(Step 9) Thus, the proposed FGP problem is constructed as follows based on model \eqref{GrindEQ__27_};
\begin{equation} \label{GrindEQ__34_} 
\begin{array}{l} {Min\, \, \beta } \\ {s.t.\left\{\begin{array}{l} {0{\rm .240}x_{1} -0{\rm .040}x_{2} +0.235x_{3} +{\rm 0.277.}+d_{1}^{-} =1,} \\ {-0.274x_{1} -0.360x_{2} -0.411x_{3} +3.074+d_{2}^{-} =1,} \\ {0\le x_{1} \le 1.051,} \\ {3.209\le x_{2} \le 8.053,} \\ {0\le x_{3} \le 1.756,} \\ {\beta \ge d_{1}^{-} ,\beta \ge d_{2}^{-} ,} \\ {\beta \ge 0,\, \, \beta \le 1,} \\ {d_{1}^{-} ,\, \, d_{2}^{-} \ge 0.} \end{array}\right. } \end{array} 
\end{equation} 
The above problem is solved and the results are obtained as:
\[x_{1} =1.051,\, \, \, x_{2} =3.331,\, \, \, x_{3} =1.432,\, \] 
\[\beta =0,\, \, d_{1}^{-} =0,\, \, d_{2}^{-} =0,\, \] 
$f_{1} =270.366$ and $f_{2} ={\rm 20.820}$ 

\noindent Also, the membership values achieved are as follows:

\noindent $\mu _{1} \left(f_{1} \right)=0.931$and$\, \mu _{2} \left(f_{2} \right)=1.00$

(Step 10): Let the decision maker is not satisfied by this solution and desires more production and then go to Step 12. According to Step 12, the new lower tolerance limit in the first objective will be 270.366. Additionally, the upper tolerance limit in the second objective function will remain the same and then from here, return to Step 5. 

\noindent Thus, membership function of the first objective is reformulated as follows (Step 6):
\[\mu _{1} \left(f_{1} \left(x\right)\right)\cong \left\{\begin{array}{l} {1\, \, \, \, \, \, \, \, \, \, \, \, \, \, \, \, \, \, \, \, \, \, \, \, \, \, \, \, \, \, \, \, \, \, \, \, \, \, \, \, \, \, \, \, f_{1} \left(x\right)\ge 281.523,} \\ {\frac{f_{1} \left(x\right)-270.366}{\left(281.523-270.366\right)} ,\, \, 270.366\le f_{1} \left(x\right)\le 281.523,\, \, \, \, \, \, \, } \\ {0\, \, \, \, \, \, \, \, \, \, \, \, \, \, \, \, \, \, \, \, \, \, \, \, \, \, \, \, \, \, \, \, \, \, \, \, \, \, \, \, \, \, \, \, 270.366\ge f_{1} \left(x\right)} \end{array}\right. \] 
In addition, membership function of the second objective function will remain unchanged.

\noindent Using the previous solution derived  in Step 9; $\left(x_{1} =1.051,\, \, \, x_{2} =3.331,\, \, \, x_{3} =1.432\right)$, the nonlinear membership function is converted into the equivalent linear function. Then the above problem \eqref{GrindEQ__34_} is reconstructed and resolved to obtain the following candidate solution: $x_{1} =1.051,\, \, \, x_{2} =3.209,\, \, \, x_{3} =1.697,\, $$\beta =\, \, d_{1}^{-} =\, \, d_{2}^{-} =0.065,\, $ $f_{1} =279.275,$ $f_{2} ={\rm 21.471,}$$\mu _{1} \left(f_{1} \right)=0.986$$\, \mu _{2} \left(f_{2} \right)=0.935.$

\noindent Let the decision maker be satisfied by the current solution $\left(x_{1} =1.051,\, \, \, x_{2} =3.209,\, \, \, x_{3} =1.697\right).$ Then the suggested algorithm is stopped at Step 11.

\noindent However, the suggested approach in [36] gives the following results: $f_{1} ={\rm 263.617}$ and$f_{2} ={\rm 22.876}$ with the membership functions achieved are as follows: $\mu _{1} \left(f_{1} \right)=0.944,$$\, \mu _{2} \left(f_{2} \right)=0.795.$ The FGP model I in [40] provides the following results:$f_{1} ={\rm 257.515,}$$f_{2} ={\rm 21.355;}$$\mu _{1} \left(f_{1} \right)=0.851,$$\, \mu _{2} \left(f_{2} \right)=0.947.$ The suggested FGP model II in [40] delivered  $f_{1} ={\rm 280.011,}$$f_{2} ={\rm 21.683;}$$\mu _{1} \left(f_{1} \right)=0.990,$$\, \mu _{2} \left(f_{2} \right)=0.914.$

\noindent From the results of the above examples, all of the sums of the membership values generated by the suggested procedure are greater than that generated by the approaches in [36] and [40]. Besides it is clear from results that the proposed procedure gives an efficient solution for the MONLPP with IT2 FNs. Moreover, it is very effective and more practical than the FGP models of papers [36] and [36] at achieving the optimal solution for the considered problems.

\section{ Conclusions}

In this paper, a type of Multiobjective Nonlinear Programming Problems (MONLPP) with trapezoidal Interval Type 2 Fuzzy Numbers (IT2 FNs) is modeled. The most serious case of the modeled problem is that it is having two objective functions; one is to maximize the desired profit while the other is to minimize the integrated time. At first, MONLPP with IT2 FNs is converted into an equivalent crisp MONLPP using an expected value function and then feasible region for the model is transformed to inequalities using the limits of decision variables. Besides membership function associated with the nonlinear objective is converted into an equivalent linear function. Thus, an interactive fuzzy goal programming approach based on Taylor series is proposed for solving the problem.

Consequently, application of the proposed procedure is discussed with two numerical models, and the effectiveness of the solutions achieved by the proposed procedure is verified. Moreover, from the above results, the suggested procedure provides an efficient solution comparing to the approaches of Mohamed [36] and Gupta et al. [40].
In future, the procedure can be extended for the multiobjective nonlinear time minimization problems. It can also be
studied for other types of uncertain programming problems like {[44],[45]} where both the membership as well as nonmembership
functions are taken into consideration.

\section{References}

\noindent \begin{enumerate}
\item [1] Zadeh, L. A. (1965). Fuzzy sets.~\textit{Information and control},~\textit{8}\eqref{GrindEQ__3_}, 338-353.

\noindent \item [2] Bellman, R. E., \& Zadeh, L. A. (1970). Decision-making in a fuzzy environment.~\textit{Management science},~\textit{17}\eqref{GrindEQ__4_}, B-141.

\noindent \item [3] Tanaka$\dagger$, H., Okuda, T., \& Asai, K. (1973). On fuzzy-mathematical programming, doi: 10.1080/01969727308545912.

\noindent \item [4] Zimmermann, H. J. (1978). Fuzzy programming and linear programming with several objective functions.~\textit{Fuzzy sets and systems},~\textit{1}\eqref{GrindEQ__1_}, 45-55.

\noindent \item [5] Nehi, H. M., \& Hajmohamadi, H. (2012). A ranking function method for solving fuzzy multi-objective linear programming problem.~\textit{Annals of Fuzzy Mathematics and Informatics},~\textit{3}\eqref{GrindEQ__1_}, 31-38.

\noindent \item [6] de Campos Ib\'{a}\~{n}ez, L. M., \& Mu\~{n}oz, A. G. (1989). A subjective approach for ranking fuzzy numbers.~\textit{Fuzzy sets and systems},~\textit{29}\eqref{GrindEQ__2_}, 145-153.

\noindent \item [7] Zadeh, L. A. (1975). The concept of a linguistic variable and its application to approximate reasoning---I.~\textit{Information sciences},~\textit{8}\eqref{GrindEQ__3_}, 199-249..

\noindent \item [8] Mendel, J. M., John, R. I., \& Liu, F. (2006). Interval type-2 fuzzy logic systems made simple.~\textit{Fuzzy Systems, IEEE Transactions on},~\textit{14}\eqref{GrindEQ__6_}, 808-821.

\noindent \item [9] Mendel, J. M. (2009). On answering the question ``Where do I start in order to solve a new problem involving interval type-2 fuzzy sets?''.~\textit{Information Sciences},~\textit{179}\eqref{GrindEQ__19_}, 3418-3431.

\noindent \item [10] Mitchell, H. B. (2005). Pattern recognition using type-II fuzzy sets.~\textit{Information Sciences},~\textit{170}\eqref{GrindEQ__2_}, 409-418.

\noindent \item [11] Zeng, W., \& Li, H. (2006). Relationship between similarity measure and entropy of interval valued fuzzy sets.~\textit{Fuzzy Sets and Systems},~\textit{157}\eqref{GrindEQ__11_}, 1477-1484.

\noindent \item [12] Wu, D., \& Mendel, J. M. (2008). A vector similarity measure for linguistic approximation: Interval type-2 and type-1 fuzzy sets.~\textit{Information Sciences},\textit{178}\eqref{GrindEQ__2_}, 381-402.

\noindent \item [13] Linda, O., \& Manic, M. (2011). Interval type-2 fuzzy voter design for fault tolerant systems.~\textit{Information Sciences},~\textit{181}\eqref{GrindEQ__14_}, 2933-2950.

\noindent \item [14] Shu, H., Liang, Q., \& Gao, J. (2008). Wireless sensor network lifetime analysis using interval type-2 fuzzy logic systems.~\textit{Fuzzy Systems, IEEE Transactions on},~\textit{16}\eqref{GrindEQ__2_}, 416-427.

\noindent \item [15] Wu, D., \& Mendel, J. M. (2007). Aggregation using the linguistic weighted average and interval type-2 fuzzy sets. Fuzzy Systems, IEEE Transactions on, 15\eqref{GrindEQ__6_}, 1145-1161.

\noindent \item [16] Han, S., \& Mendel, J. M. (2010, July). Evaluating location choices using perceptual computer approach. In~\textit{Fuzzy Systems (FUZZY), 2010 IEEE International Conference on}~(pp. 1-8). IEEE.

\noindent \item [17] Chen, S. M., \& Lee, L. W. (2010). Fuzzy multiple attributes group decision-making based on the ranking values and the arithmetic operations of interval type-2 fuzzy sets.~\textit{Expert Systems with applications},~\textit{37}\eqref{GrindEQ__1_}, 824-833.

\noindent \item [18] Sinha, B., Das, A., \& Bera, U. K. (2015). Profit maximization solid transportation problem with trapezoidal interval type-2 fuzzy numbers. \textit{International Journal of Applied and Computational Mathematics}, 1-16.

\noindent \item [19] Li, H., Pan, Y., \& Zhou, Q. (2015). Filter design for interval type-2 fuzzy systems with D stability constraints under a unified frame.~\textit{Fuzzy Systems, IEEE Transactions on},~\textit{23}\eqref{GrindEQ__3_}, 719-725.

\noindent \item [20] Hu, J., Zhang, Y., Chen, X., \& Liu, Y. (2013). Multi-criteria decision making method based on possibility degree of interval type-2 fuzzy number. .\textit{Knowledge-Based Systems},~\textit{43}, 21-29.

\noindent \item [21] Qin, J., Liu, X., \& Pedrycz, W. (2015). An extended VIKOR method based on prospect theory for multiple attribute decision making under interval type-2 fuzzy environment.~\textit{Knowledge-Based Systems},~\textit{86}, 116-130.

\noindent \item [22] Kahraman, C., \"{O}ztay\c{s}i, B., Sar$\imath$, \.{I}. U., \& Turano\u{g}lu, E. (2014). Fuzzy analytic hierarchy process with interval type-2 fuzzy sets. Knowledge-Based Systems, 59, 48-57.

\noindent \item [23] Sari, I. U., \& Kahraman, C. (2015). Interval Type-2 Fuzzy Capital Budgeting, \textit{International Journal of Fuzzy Systems},~\textit{17}\eqref{GrindEQ__4_}, 635-646.

\noindent \item [24] Cevik Onar, S., Oztaysi, B., \& Kahraman, C. (2014). Strategic decision selection using hesitant fuzzy TOPSIS and interval type-2 fuzzy AHP: a case study.~\textit{International Journal of Computational intelligence systems},~\textit{7}\eqref{GrindEQ__5_}, 1002-1021.

\noindent \item [25] Bustince Sola, H., Fernandez, J., Hagras, H., Herrera, F., Pagola, M., \& Barrenechea, E. (2015). Interval type-2 fuzzy sets are generalization of interval-valued fuzzy sets: toward a wider view on their relationship. Fuzzy Systems, IEEE Transactions on, 23\eqref{GrindEQ__5_}, 1876-1882.

\noindent \item [26]  Castillo, O., \& Melin, P. (2014). A review on interval type-2 fuzzy logic applications in intelligent control.~\textit{Information Sciences},~\textit{279}, 615-631.

\noindent \item [27] Sanchez, M. A., Castillo, O., \& Castro, J. R. (2015). Generalized type-2 fuzzy systems for controlling a mobile robot and a performance comparison with interval type-2 and type-1 fuzzy systems.~\textit{Expert Systems with Applications},\textit{42}\eqref{GrindEQ__14_}, 5904-5914.

\noindent \item [28] Gaxiola, F., Melin, P., Valdez, F., \& Castillo, O. (2014). Interval type-2 fuzzy weight adjustment for backpropagation neural networks with application in time series prediction.~\textit{Information Sciences},~\textit{260}, 1-14.

\noindent \item [29] Soto, J., Melin, P., \& Castillo, O. (2014). Time series prediction using ensembles of ANFIS models with genetic optimization of interval type-2 and type-1 fuzzy integrators.~\textit{International Journal of Hybrid Intelligent Systems},\textit{11}\eqref{GrindEQ__3_}, 211-226.

\noindent \item [30] Pulido, M., Melin, P., \& Castillo, O. (2014). Particle swarm optimization of ensemble neural networks with fuzzy aggregation for time series prediction of the Mexican Stock Exchange.~\textit{Information Sciences},~\textit{280}, 188-204.

\noindent \item [31] Melin, P., \& Castillo, O. (2014). A review on type-2 fuzzy logic applications in clustering, classification and pattern recognition.~\textit{Applied soft computing},~\textit{21}, 568-577.

\noindent \item [32] Melin, P., Gonzalez, C. I., Castro, J. R., Mendoza, O., \& Castillo, O. (2014). Edge-detection method for image processing based on generalized type-2 fuzzy logic.~\textit{Fuzzy Systems, IEEE Transactions on},~\textit{22}\eqref{GrindEQ__6_}, 1515-1525.

\noindent \item [33] Narasimhan, R. (1980). Goal programming in a fuzzy environment.~\textit{Decision sciences},~\textit{11}\eqref{GrindEQ__2_}, 325-336.

\noindent \item [34] Hannan, E. L. (1981). On fuzzy goal programming*.~\textit{Decision Sciences},~\textit{12}\eqref{GrindEQ__3_}, 522-531.

\noindent \item [35] Tiwari, R. N., Dharmar, S., \& Rao, J. R. (1987). Fuzzy goal programming---an additive model.~\textit{Fuzzy sets and systems},~\textit{24}\eqref{GrindEQ__1_}, 27-34.

\noindent \item [36] Mohamed, R. H. (1997). The relationship between goal programming and fuzzy programming.~\textit{Fuzzy Sets and Systems},~\textit{89}\eqref{GrindEQ__2_}, 215-222.

\noindent \item [37] Pal, B. B., Moitra, B. N., \& Maulik, U. (2003). A goal programming procedure for fuzzy multiobjective linear fractional programming problem.~\textit{Fuzzy Sets and Systems},~\textit{139}\eqref{GrindEQ__2_}, 395-405.

\noindent \item [38] Toksar$\imath$, M. D., \& Bilim, Y. (2015). Interactive Fuzzy Goal Programming Based on Jacobian Matrix to Solve Decentralized Bi-level Multi-objective Fractional Programming Problems.~\textit{International Journal of Fuzzy Systems},~\textit{17}\eqref{GrindEQ__4_}, 499-508.

\noindent \item [39] Dalman, H., K\"{o}\c{c}ken, H. G., \& Sivri, M. (2013). A Solution Proposal to Indefinite Quadratic Interval Transportation Problem.~\textit{New Trends in Mathematical Sciences},~\textit{1}\eqref{GrindEQ__2_}, 07-12.

\noindent \item [40] Gupta, M., \& Bhattacharjee, D. (2012). Two weighted fuzzy goal programming methods to solve multiobjective goal programming problem.~\textit{Journal of Applied Mathematics,~2012}, doi: 10.1155/2012/796028

\noindent \item [41] Tang, J., \& Wang, D. (1996, December). Modelling and optimization for a type of fuzzy nonlinear programming problems in manufacturing systems. In \textit{Decision and Control, 1996. Proceedings of the 35th IEEE Conference on}~(Vol. 4, pp. 4401-4405). IEEE.

\noindent \item [42] Singh, S. K., \& Yadav, S. P. (2015). Intuitionistic fuzzy non linear programming problem: Modeling and optimization in manufacturing systems.~\textit{Journal of Intelligent \& Fuzzy Systems},~\textit{28}\eqref{GrindEQ__3_}, 1421-1433.

\noindent \item [43] Zangiabadi, M., \& Maleki, H. R. (2013). Fuzzy goal programming technique to solve multiobjective transportation problems with some non-linear membership functions.~\textit{Iranian Journal of Fuzzy Systems},~\textit{10}\eqref{GrindEQ__1_}, 61-74.

\noindent \item [44] Dalman, H. (2016). Uncertain programming model for multi-item solid transportation problem. International Journal of Machine Learning and Cybernetics, doi: 10.1007/s13042-016-0538-7.

\noindent \item [45] Dalman, H., Guzel, N., and Sivri, M. (2015) A Fuzzy Set-Based Approach to Multi-objective Multi-item Solid Transportation Problem Under Uncertainty. International Journal of Fuzzy Systems, doi: 10.1007/s40815-015-0081-9.
\end{enumerate}

\noindent

\end{document}